\documentclass[11pt,english,reqno]{amsart}

\usepackage[top=3cm,bottom=3.5cm,inner=3cm,outer=3cm]{geometry}
\usepackage{amsthm}
\usepackage{amsmath}
\usepackage{amssymb}
\usepackage{setspace}
\usepackage[numbers]{natbib}
\usepackage{mathtools}
\usepackage{graphics}
\usepackage{epstopdf}
\usepackage[hidelinks]{hyperref}
\usepackage{enumerate}
\usepackage{csquotes}
\usepackage{verbatim}
\usepackage{tikz}
\usetikzlibrary{calc}
\usetikzlibrary{decorations.pathreplacing}

\MakeOuterQuote{"}

\makeatletter
\def\imod#1{\allowbreak\mkern10mu({\operator@font mod}\,\,#1)}

\makeatother

\theoremstyle{plain}
\newtheorem*{thm*}{Theorem}
\newtheorem{thm}{Theorem}
\newtheorem{lem}[thm]{Lemma}
\newtheorem{claim}[thm]{Claim}

\newtheorem*{claim*}{Claim}
\newtheorem{cor}[thm]{Corollary}

\theoremstyle{definition}

\theoremstyle{remark}

\usepackage{babel}

\def\eps{\varepsilon}

\begin{document}

\title{Path Ramsey number for random graphs}

\author[S. Letzter]{Shoham Letzter}
\address{Department of Pure Mathematics and Mathematical Statistics, University of Cambridge, Wilberforce Road, Cambridge CB3\thinspace0WB, UK}
\email{s.letzter@dpmms.cam.ac.uk}

\date{}

\begin{abstract}
\setlength{\parskip}{\medskipamount}
\setlength{\parindent}{0pt}
\setstretch{1.2}

Answering a question raised by Dudek and Pra\l{}at \cite{pralat}, we show that if $pn\rightarrow \infty$, w.h.p.,~whenever $G=G(n,p)$ is $2$-coloured, there exists a monochromatic path of length $n(2/3+o(1))$. This result is optimal in the sense that  $2/3$ cannot be replaced by a larger constant.

As part of the proof we obtain the following result which may be of independent interest.
We show that given a graph $G$ on $n$ vertices with at least $(1-\eps)\binom{n}{2}$ edges, whenever $G$ is $2$-edge-coloured, there is a monochromatic path of length at least $(2/3-100\sqrt{\eps})n$. This is an extension of the classical result by Gerencs\'er and Gy\'arf\'as \citep{gyarfas} which says that whenever $K_n$ is $2$-coloured there is a monochromatic path of length at least $2n/3$.

\end{abstract}

\maketitle

\section{Introduction}
Considering the richness of Ramsey theory and the great interest in random graphs, it is natural to consider Ramsey properties of random graphs.

The study of random Ramsey theory has proved particularly useful in the establishment of upper bounds on the size Ramsey number.
For graphs $G,F,H$, we write $G\rightarrow (F,H)$ if for every red-blue colouring of the edges of $G$, there is either a red $F$ or a blue $H$. If $F,H$ are isomorphic, we use instead the notation $G\rightarrow H$. The size Ramsey number, denoted by $\hat{r}(H)$ is defined to be $\hat{r}(H)=\min\{|E(G)|: G\rightarrow H\}$.

In \cite{beck}, disproving a conjecture of Erd\H{o}s \cite{erdos}, Beck showed that $\hat{r}(P_n)\le 900 n$.  In \citep{bela} Bollob\'as noted a slightly better bound, and recently Dudek and  Pra\l{}at \citep{pralat} gave an elementary proof of the bound $\hat{r}(P_n)\le 137n$. They actually proved that w.h.p.,~$G(n,\alpha/n)\rightarrow P_{\beta n}$ for some constants $\alpha, \beta$.
They raised the question of determining the maximum $l$ such that $G(n,p)\rightarrow P_l$, where $pn\rightarrow \infty$.
Inspired by the well known result of  Gerencs\'er and Gy\'arf\'as \cite{gyarfas} which says that $K_n\rightarrow P_{2n/3}$, they ask if $G(n,p)\rightarrow P_l$ for some $l=n(2/3+o(1))$.
Our main result answers this question in the affirmative.

\begin{thm}\label{thm_random_path_ramsey}
Let $0<p=p(n)<1$ and assume that $pn\rightarrow\infty$.
Then w.h.p.,  $G(n,p)\rightarrow P_l$ for some  $l=(2/3+o(1))n$.
\end{thm}

This result is essentially best possible since there is a $2$-colouring of the edges of $K_n$ such that the longest monochromatic path is of length $\lceil 2n/3 +1\rceil$. To see this, divide the vertex set of $K_n$ into two sets $A,B$ such that $|A|=\lfloor n/3 \rfloor$, let the edges spanned by $B$ be coloured red and colour the other edges blue.

In order to prove Theorem \ref{thm_random_path_ramsey}, we prove the following extension of the result in \cite{gyarfas} to graphs with a large number of edges.

\begin{thm}\label{thm_density_path_ramsey}
Let $0\le\eps\le 1/4$, $k\ge l$ and let $G$ be a graph on $n\ge k+\lfloor(l+1)/2\rfloor+150\sqrt{\eps} k$ vertices, with at least $(1-\eps)\binom{n}{2}$ edges.
Then $G\rightarrow (P_{k+1},P_{l+1})$.

In particular, given $0\le \eps<1/90000$,  for every graph $G$ on $n$ vertices and at least $(1-\eps)\binom{n}{2}$ edges,  $G\rightarrow P_k$ where $k=(2/3-100\sqrt{\eps})n$.

\end{thm}

Theorem \ref{thm_density_path_ramsey} is a consequence of the following similar result, in which we consider graphs with large minimum degree rather than high density.

\begin{thm}\label{thm_min_deg_path_ramsey}
Let $0<\eps<1/4$, $k\ge l$ and let $G$ be a graph on $n\ge k+\lfloor(l+1)/2\rfloor+60\eps k$  vertices with minimum degree at least $(1-\eps)n$. Then  $G\rightarrow (P_{k+1},P_{l+1})$.
\end{thm}

Note that it is easy to deduce Theorem \ref{thm_density_path_ramsey} from Theorem \ref{thm_min_deg_path_ramsey}.
By an averaging argument, it suffices to prove the assertion for $n=k+l/2+150\sqrt{\eps}$.
By removing at most $\sqrt{\eps}n$ vertices, we obtain a graph on $n'\ge (1-\sqrt{\eps})n$ vertices and minimum degree at least $(1-\sqrt{\eps})n\ge (1-\sqrt{\eps})n'$ vertices.
One can check that $(1-\sqrt{\eps})n\ge k+l/2+100\sqrt{\eps}k$, so the assertion of Theorem \ref{thm_density_path_ramsey} follows from Theorem \ref{thm_min_deg_path_ramsey}.
For the second part, it is easy to check that when $k=l=(2/3-100\sqrt{\eps})n$ and $150\sqrt{\eps}\le 1/2$, it follows that $n\ge k+l/2+150\sqrt{\eps}$.

The rest of the paper is organised as follows.
In Section \ref{sec_dense} we prove Theorem \ref{thm_min_deg_path_ramsey}. In 
order to prove Theorem \ref{thm_random_path_ramsey}, we use the so-called sparse regularity lemma, due to Kohayakawa \citep{kohayakawa} and R\"odl (see \cite{conlon}). In Section \ref{sec_regularity} we state this result as well as some necessary notation. 
We prove Theorem \ref{thm_random_path_ramsey} in Section \ref{sec_random} and finish with some concluding remarks in Section \ref{sec_conclusion}.
Throughout the paper we omit floor and ceiling signs  whenever they do not affect the arguments.

\section{Path Ramsey number for dense graphs}\label{sec_dense}
In the proof of Theorems \ref{thm_min_deg_path_ramsey} and \ref{thm_random_path_ramsey} we use the following  observation, from \cite{pralat} and \citep{Pokrovskiy}. For the sake of completeness, we prove it here.
\begin{lem}\label{lem_partition}
For every graph $G$ there exist two disjoint sets of vertices $U,W$ of equal sizes, such that there are no edges between them and $G\backslash (U\cup W)$ has a Hamiltonian path.
\end{lem}

\begin{proof}
In order to find sets with the desired properties, we apply the following algorithm, maintaining a partition of $V(G)$ into sets $U,W$ and a path $P$.
Start with $U=V(G), W=\emptyset$ and $P$ an empty path.
At every stage in the algorithm,  do the following.
If $|U|\le|W|$,  stop.
Otherwise, if $P$ is empty,  move a vertex from $U$ into $W$ (note that $U\neq \emptyset$).
If $P$ is non-empty, let $v$ be its endpoint. If $v$ has a neighbour $u$ in $U$,  put $u$ in $P$, otherwise  move $v$ to $W$.

Note that at any given point in the algorithm there are no edges between $U$ and $W$.
Furthermore, the value $|U|-|W|$ is positive at the beginning of the algorithm and decreases by one at every stage, thus at some point the algorithm will stop and produce sets $U,W$ with the required properties.
\end{proof}

Occasionally it is easier to use the following immediate consequence of Lemma \ref{lem_partition}.

\begin{cor}\label{cor_partition}
Let $G$ be a balanced bipartite graph on $n$ vertices with bipartition $V_1,V_2$, which has no path of length $k$.  Then there exist $X_i\subseteq V_i$ such that $|X_1|=|X_2|\ge (n-k)/4$ and $G$ has no edges between $X_1$ and $X_2$.
\end{cor}

\begin{proof}
Let $U,W$ be as in Lemma \ref{lem_partition} and let $P$ be a Hamiltonian path in $G\backslash(U\cup W)$.
Note that $P$ must alternate between $V_1$ and $V_2$, thus $|V(P)\cap V_1|=|V(P)\cap V_2|$ (it follows from the assumptions that the number of vertices in $P$ is even).
Denote $U_i=U\cap V_i$, $W_i=W\cap V_i$, for $i=1,2$ and assume that $|U_1|\ge |U_2|$.
It follows that $|U_1|+|W_1|=|U_2|+|W_2|$.
Thus, using the fact that $|U|=|W|$, we have that $|U_1|=|W_2|\ge |U|/2\ge (n-k)/4$.
Set $X_1=U_1$ and $X_2=W_2$.
\end{proof}

We now prove Theorem \ref{thm_min_deg_path_ramsey}.

\begin{proof} [Proof of Theorem \ref{thm_min_deg_path_ramsey}]

We prove the theorem by induction on $k$.
Clearly, if $k=1$ the claim holds.
Given $k>1$, let $G$ be a graph on $n\ge  k+\lfloor (l+1)/2\rfloor+100\eps k$ vertices, with  minimum degree at least $(1-\eps)n$, and consider a red-blue colouring of the edges of $G$.

If $k>l$ then by induction there is either a red $P_k$ or a blue $P_{l+1}$; in the latter case we are done.
If $k=l$ then by induction there is either a red or blue $P_k$.
Thus, without loss of generality there is a red path of length $k-1$, which we denote by $P=(v_1,\ldots,v_k)$. Let $U=V(G)\backslash V(P)$.

We note first that the assertion of Theorem \ref{thm_min_deg_path_ramsey} holds when $k\ge n(1/2-\eps)$.
If there is no red $P_{k+1}$, then by Lemma \ref{lem_partition}, we can find disjoint sets $U,W$, of size at least $(n-k)/2$ such that there is no red edge between them.
Since $G$ has minimum degree at least $n(1-\eps)$, we can greedily find a blue path of length at least $|U|+|W|-2\eps n\ge n(1-2\eps)-k\ge k$.
Thus we can assume that $k\ge 4n$, so every vertex in $G$ has at most $4\eps k$ non neighbours.
Put $\delta=4\eps$.
Note that we can assume that $\delta k\ge 1$, otherwise $G$ is a complete graph and Theorem \ref{thm_min_deg_path_ramsey} follows directly from \cite{gyarfas}.

We consider three  cases.

\subsection{Case 1. $G[U]$ contains a blue path $Q$ of length $13\delta k$}

Let $Q_1$ be a maximal path extending $Q$ by alternating between vertices of $P$ and $U$ and which has both ends in $U$.
Let $U'=U\backslash V(Q_1)$ and $V'=V(P)\backslash V(Q_1)$.
Let $Q_2$ be a maximal path alternating between $U'$ and $V'$ which has both ends in $U$.
Denote the ends of $Q_i$ by $x_i,y_i$, for $i=1,2$.
We show that $|Q_1|+|Q_2|\ge l+3\delta k$.

Suppose this is not the case. In particular,  $Q_1,Q_2$ do not cover $U$, so we can pick a  vertex $z\in U\backslash (V(Q_1)\cup V(Q_2) )$.
Note that all but at most $3\delta k$ vertices of $P$ are adjacent to all of $x_1,x_2,z$. By our assumption on the lengths of $Q_1$ and $Q_2$, the number of vertices of $P$ which are included in one of $Q_1$ and $Q_2$ is at most $k/2 -5\delta k$, hence there exist vertices $v_i,v_{i+1}$ which are adjacent to all of $x_1,y_1,z$.
We assume that $v_i$ and $v_{i+1}$ have no common red neighbour in $x_1,x_2,z$ because otherwise we obtain a red $P_{k+1}$.
It follows that without loss of generality, $v_i$ is joined in blue to two of $x_1,y_1,z$, contradicting the maximality of $Q_1$ and $Q_2$.

Let $Q_2'$ be a subpath of $Q_2$ with ends $x_2',y_2'\in U$ satisfying $|Q_2|+|Q_1|=l+3\delta k$.
A similar argument to the above shows that without loss of generality there exist $v_i,v_{i+1}$ such that $x_1,y_1$ are blue neighbours of $v_i$ and $x_2',y_2'$ are blue neighbours of $v_{i+1}$.
Denote by $C_1$ and $C_2$ the blue cycles obtained by adding $v_i$ to $Q_1$ and $v_{i+1}$ to $Q_2$, and let $U_i=V(C_i)\cap U$.

Note that if there is any blue edge between $C_1$ and $C_2$ we obtain a blue path of length $l$, so we assume that no such edges exist. We can also assume that $|U_1|,|U_2|\ge 3\delta k$, otherwise one of $Q_1,Q_2$ has length at least $l$.

The number of vertices in $V(P)\backslash (V(C_1)\cup V(C_2))$ is at least $k/2+5\delta k$, hence there exists $j$ such that $v_j,v_{j+1}\notin V(C_1)\cup V(C_2)$.
If one of $v_j$ and $v_{j+1}$ has blue neighbours in both $U_1$ and $U_2$, we obtain a blue path of length $l$, so we can assume this is not the case.
Also, we assume that $v_j,v_{j+1}$ have no red common neighbour in either $U_1$ or $U_2$, because otherwise we obtain a red path of length $k$.
Thus, recalling that $v_j,v_{j+1}$ have at most $\delta k$ non neighbours in $G$, without loss of generality, $v_j$ is joined in red to all but $\delta k$ vertices of $U_1$, and $v_{j+1}$ is joined in red to all but $\delta k$ vertices in $U_2$.
Let $w_1\in U_1$ be any red neighbour of $v_j$. Since it is connected to all but at most $\delta k$ vertices of $U_2$ and these edges must all be red, $U_2$ contains a vertex $w_2$ which is a red neighbour of both $w_1$ and $v_{j+1}$.
We obtain a red path $v_1,\ldots,v_j,w_1,w_2,v_{j+1},\ldots,v_k$ of length $k$.

This finishes the proof of Theorem \ref{thm_density_path_ramsey} in the first case.

\subsection{Case 2. $l\le (1-13\delta) k$}

Let $Q_1$ be a maximal blue path alternating between $U$ and $P$ and having both ends in $U$ and similarly  let $Q_2$ be a maximal blue path alternating between $U\backslash V(Q_1)$ and $V(P)\backslash V(Q_1)$.
As in the previous case, it can be shown that $|Q_1|+|Q_2|\ge l+3\delta k$.

Let $Q_2'$ be a subpath of $Q_2$ such that $|Q_1|+|Q_2'|=l+3\delta k$.
As before, there exists $j$ such that both $v_j,v_{j+1}\in V(P)\backslash (V(Q_1)\cup V(Q_2'))$ and they are joined in $G$ to all ends of the two paths. Thus the vertices $v_j,v_{j+1}$ can be used to extend $Q_1,Q_2'$ into blue vertex disjoint cycles $C_1,C_2$, whose sum of length is $l+3\delta k$ and each of which has length at least $3\delta k$.
The proof of Theorem \ref{thm_min_deg_path_ramsey} can now be finished as in the first case. 

\subsection{Case 3. $l\ge (1-13\delta)k$ and $G[U]$ contains no blue path of length at least $13\delta k$}

We conclude from Lemma \ref{lem_partition} that there exist two disjoint sets $W_1,W_2\subseteq U$ of size $|W_1|=|W_2|\ge (1/2+3\delta) k/2$ with no blue edges between them. 
Since every vertex in $G$ is adjacent to all but at most $\delta k$ vertices, we can greedily find a red path $Q$ in $U$ of length at least $|W_1|+|W_2|-2\delta k=(1/2 +\delta )k$.

Let $X$ be the set of the first and last $ (1/4+\delta/2)k$ vertices of $P$.
We  assume that there is no red edge between $X$ and $Q$, because otherwise there is a red path of length $k$.
We can now greedily construct a blue path alternating between $X$ and $V(Q)$ of length at least $|X|+|Q|-2\delta k\ge k\ge l$.
\end{proof}

\section{Sparse regularity lemma}\label{sec_regularity}

We shall make use of a variant of Szemer\'edi's regularity lemma  \cite{regularity} for sparse graphs, often referred to as the sparse regularity lemma, which was proved independently by Kohayakawa \cite{kohayakawa} and R\"odl (see \cite{conlon}).
Before stating the theorem, we introduce some notation.

Given two disjoint sets of vertices $U,V$ in a graph, we define the density $d_p(U,V)$ of edges between $U$ and $V$ with respect to $p$ to be
\begin{equation}
d_p(U,V)=\frac{e(U,V)}{p|U||V|},
\end{equation}
 where $e(U,V)$ is the number of edges between $U$ and $V$.
We say that a bipartite graph with bipartition $U,V$ is $(\eps, p)$-regular if for every $U'\subseteq U, V'\subseteq V$ with $|U'|\ge \eps |U|, |V'|\ge \eps|V| $ the density $d_p(U',V')$ satisfies $|d_p(U',V')-d_p(U,V)|\le \eps $.

Given a graph $G$, a partition $V_1,\ldots,V_t$ of $V(G)$ is called an $(\eps,p)$-regular partition if it is an equipartition (i.e. the sizes of the sets differ by at most one), and if all but at most $\eps$ of the pairs $V_i,V_j$ induce an $(\eps,p)$-regular graph.

Given $0<\eta,p<1, D\ge 1$, a graph $G$ is called $(\eta, p, D)$-upper-uniform if for all disjoints sets of vertices $U_1,U_2$ of size at least $\eta |V(G)|$, the density $d_p(U_1,U_2)$ is at most $D$. Note that random graphs are w.h.p.~upper uniform (with suitable parameters).

We are now ready to state the sparse regularity lemma of Kohayakawa and R\"odl.
\begin{thm}\label{thm_sparse_regularity}
For every $\eps>0$, $t$ and $D>1$ there exist $\eta>0$ and $T$ such that for every $0\le p\le 1$, every $(\eta,p,D)$-upper-uniform graph admits an $(\eps,p)$-regular partition into $s$ parts where $t\le s\le T$.
\end{thm}

We shall use a slightly stronger variant of \ref{thm_sparse_regularity}, namely the coloured version of the sparse regularity lemma.
\begin{thm}\label{thm_coloured_sparse_regularity}
For every $\eps>0$, $t,l$ and $D>1$ there exist $\eta>0$ and $T$ such that for every $0\le p\le 1$, if $G_1,\ldots,G_l$ are $(\eta,p,D)$-upper-uniform graphs on vertex set $V$, there is an equipartition of $V$ into $s$ parts, where $t\le s\le T$, for which all but at most $\eps$ of the pairs induce a regular pair in each $G_i$.
\end{thm}

\section{Path Ramsey number for random graphs}\label{sec_random}
Before turning to the proof of Theorem \ref{thm_random_path_ramsey}, we remark that a weaker result can be proved using elementary tools.
\begin{lem}
Let $0<p=p(n)<1$ and assume that $pn\rightarrow\infty$.
Then w.h.p.,  $G(n,p)\rightarrow P_l$ for some  $l=(1/2+o(1))n$.
\end{lem}

\begin{proof}
Given $\alpha>0$, suppose $G$ can be coloured such that there is no monochromatic path of length $n(1/2-\alpha)$.
By Lemma \ref{lem_partition}, there exist disjoint sets $U,W$, both of size at least $n(1/2+\alpha)/2$ with no red edges between them.
Considering the graph $G[U,W]$, it follows from the same lemma that there exist disjoint sets $X\subseteq U,Y\subseteq W$ of size at least $\alpha n/2$, such that there are no blue edges between them. We conclude that there are no edges of $G$ between $X$ and $Y$.
But w.h.p.,~every two disjoint sets of at least $\alpha n/2$ vertices in $G$ have an edge between them.
This shows that w.h.p., in every $2$-colouring of $G$ there is a monochromatic path of length at least $n(1/2-\alpha)$.
\end{proof}

\begin{proof}[Proof of Theorem \ref{thm_random_path_ramsey}]
Let $0<p<1$ be such that $pn\rightarrow \infty$ and let $\alpha>0$.

We show that w.h.p., for every $2$-edge-colouring of $G=G(n,p)$ there is a monochromatic path of length at least $(2/3-\alpha)n$.
We pick $\eps>0$ small and $t$ large (taking $t=1/\eps$ and $\eps$ small enough such that $60\sqrt{\eps}\le \alpha$ would do).

Let $\eta, T$ be the constants arising from the application of Theorem \ref{thm_coloured_sparse_regularity} with $\eps,t,l=2,D=2$.
Note that $G$ is w.h.p.~$(\eta,p,2)$-upper-uniform. Thus, by Theorem \ref{thm_coloured_sparse_regularity}, given a $2$-edge-colouring of $G$, there exists an $(\eps,p)$-regular partition $V_1,\ldots,V_s$  with $t\le s\le T$.
Again, w.h.p.,~the density of edges $d_p(V_i,V_j)$ is at least $1/2$.

Let $H$ be the auxiliary graph with vertex set $[s]$ where $ij$ is an edge iff $V_i,V_j$ induce a regular bipartite graph in both red and blue.
We colour an edge $ij$ in $H$ red if the red density $d_p(V_i,V_j)$ is at least $1/4$ and blue otherwise, so if $ij$ is blue then the blue density is at least $1/4$.

Since the partition  $V_1,\ldots,V_s$ is $(\eps,p)$-regular, the number of edges in $H$ is at least $(1-\eps)\binom{s}{2}$.
It follows from Theorem \ref{thm_density_path_ramsey} that $H$ contains a monochromatic path $P$ on at least $l=(2/3-\delta)s$ vertices, where $\delta=50\sqrt{\eps}$ (assuming $\eps>0$ is small enough).
Denote by $i_1,\ldots,i_l$ the vertices of $P$.

Assuming without loss of generality that $P$ is red, we show that  $G$ contains a red path of length at least $(2/3-\alpha)n$.
We divide each set $V_{i_j}$ into two sets $U_j,W_j$ of equal sizes, so $|U_j|=n/2s$.
Let $P_j$ be a longest red path in the bipartite graph $G[U_j,W_{j+1}]$. In the following claim we show that $P_j$ covers most vertices in $U_j\cup W_{j+1}$. We shall then show that consecutive paths $P_j,P_{j+1}$ can be connected without losing too many vertices, thus obtaining a red path in $G$ of the required length.

\begin{claim}\label{claim_path_regular_pair}
For every $1\le j\le l$, $P_j$ covers at least $1-4\eps$ of the vertices of $U_j\cup W_{j+1}$.
\end{claim}

\begin{proof}
Suppose that for some $j$, $P_j$ covers at most $1-4\eps$ of the vertices of $U_j\cup W_{j+1}$. 
Set $U=U_j$ and $W=W_{j+1}$.
By Corollary \ref{cor_partition}, there exist  sets $X\subseteq U, Y\subseteq W$ with $|X|=|Y|\ge \eps|U|$, such that there are no red edges between $X$ and $Y$.
But by the regularity of the partition $V_1,\ldots,V_s$, the density $d_p(U,V)$ is within $\eps $ of the density of red edges between $U$ and $W$, which is at least $1/4$. In particular, $G$ has a red edge between $X$ and $Y$, contradicting our assumption, so Claim \ref{claim_path_regular_pair} holds.
\end{proof}

We now show that the paths $P_1,\ldots, P_{l-1}$ can be joined to a path $Q$ without losing many of the vertices.
Let $X_j$ be the set of first $2\eps |V_1|$ vertices of $P_j$ and similarly let $Y_j$ be the set of last $2\eps |V_1|$ vertices of $P_j$.
Since the paths $P_j$ alternate between the sets $U_j,W_{j+1}$, we have that $|Y_j\cap V_{i_j}|,|X_{j+1}\cap V_{i_{j+1}}|\ge \eps |V_1|$. It follows from the fact that $i_ji_{j+1}$ is a blue edge in $H$ that there is a blue edge between $Y_j$ and $X_{j+1}$.
Hence $G$ has a blue path $Q$ which contains all vertices of $V(P_1)\cup\ldots\cup V(P_{l-1})$ but at most $4\eps |V_1|(l-1)$. 

Using Claim \ref{claim_path_regular_pair}, we have that 
$|P_j|\ge (1-4\eps)|V_1|$, so 
\begin{align*}
|Q|\ge &(1-8\eps)(l-1)|V_1|=
(1-8\eps)(s(2/3-\delta)-1)n/s\ge\\
& n(2/3-(\delta+1/t+6\eps))\ge n(2/3-\alpha).
\end{align*}

This completes the proof of Theorem \ref{thm_random_path_ramsey}.

\end{proof}

\section{Concluding Remarks}\label{sec_conclusion}
It is easy to construct examples of graphs $G$ on $n$ vertices with $n\ge k+\lfloor(l+1)/2\rfloor+c\eps k$ and at least $(1-\eps)\binom{n}{2}$ edges which admits a red-blue colouring with no red $P_{k+1}$ or blue $P_{l+1}$  (e.g. by letting as many vertices as possible be isolated).
It may be interesting to determine the correct dependence of $n$ on $\eps$ in Theorem \ref{thm_density_path_ramsey}. In particular, can $\sqrt{\eps}$  be replaced by a factor of $\eps$?
Similarly, it would be interesting to determine if the error term $c\eps k$ in Theorem \ref{thm_min_deg_path_ramsey} can be replaced by $o(\eps k)$.

\bibliography{path}
\bibliographystyle{siam}

\end{document}